\documentclass[11pt]{amsart}
\usepackage{geometry}                % See geometry.pdf to learn the layout options. There are lots.
\usepackage{graphicx}
\usepackage{amssymb}
\usepackage{epstopdf}
\usepackage{color}
\title[The counterexample which I missed]{A note on a recent result by Lorist, Meyries, and Veraar, or\\ the counterexample which I missed.}
\author{Hans Zwart}
\date{September, 24 -- 28, 2026. Section 4 was added later.}                                           % Activate to display a given date or no date

\begin{document}
\begin{abstract}
In \cite{LoMV26} Lorist, Meyries, and Veraar show, among others, that the inverse of an infinitesimal generator of a bounded semigroup on a Hilbert space need not to be an infinitesimal generator of a bounded semigroup. In this note we show that with the ingredients of that paper, the example can be constructed more directly. Furthermore, the author reflects on the reason why he did not construct this a decade ago.
\end{abstract}
\maketitle

\section{Introduction}

The inverse generator problem concerns the question whether the inverse of the infinitesimal generator of a bounded semigroup on a Banach space $X$ is itself the generator of a bounded semigroup. The problem was first posed by deLaubenfels in \cite{deLa88}. It is easy to see that this assertion holds when the state space $X$ is finite-dimensional. In \cite{GoZT07} it is shown that the assertion does not hold in Banach spaces, even not for $X= \ell^p$, $p \neq 2$. The case when $X$ is a Hilbert space remained open for more than a decade. Recently, Lorist, Meyries, and Veraar \cite{LoMV26}  showed that the assertion does not hold in Hilbert spaces either. Based on their construction, a more direct example will be constructed in the next section. 

\section{Counter example}

As in \cite{LoMV26}, we define the sequences
\begin{alignat}{2}
\label{eq:y}
  y_1 &=1, \quad &&y_{k+1} = y_k^2 +1,\\
%\end{equation}
%\begin{equation}
\label{eq:lambda}
  \mu_k &= -1 + i y_k , \qquad &&\lambda_k = \mu_k^{-1}.
\end{alignat}
In the proof of their Lemma 2.2, equation (2.5), it is shown that 
\begin{equation}
\label{eq:var}
  \sum_{k=1}^{\infty} \left| e^{\lambda_{k+1} t} - e^{\lambda_{k} t}\right| \leq 8, \quad \forall t\geq 0.
\end{equation}
Let $\{\phi_k, k \in {\mathbb N}\}$ be a conditional basis of the Hilbert space $X$, which is not unconditional. We define for this basis the (linear) operator
\begin{equation}
\label{eq:eAt}
  T(t)\sum_{k=1}^{k_e} \alpha_k \phi_k = \sum_{k=1}^{k_e} \alpha_k e^{\lambda_k t}\phi_k, \quad k_e \in {\mathbb N} .
\end{equation}
By (\ref{eq:var}) it follows that $T(t)$ extends to a bounded operator on $X$ which is uniformly bounded (in $t\geq 0$). Furthermore, it is a strongly continuous semigroup with the bounded generator, $A$ given by
\begin{equation}
\label{eq:A}
  A\sum_{k=1}^{k_e} \alpha_k \phi_k = \sum_{k=1}^{k_e} \alpha_k \lambda_k \phi_k, \quad k_e \in {\mathbb N} .
\end{equation}
It is easy the see that its inverse exists as a densely defined (unbounded) operator and satisfies
\begin{equation}
\label{eq:Ainv}
  A^{-1}\sum_{k=1}^{k_e} \alpha_k \phi_k = \sum_{k=1}^{k_e} \alpha_k \lambda_k^{-1} \phi_k = \sum_{k=1}^{k_e} \alpha_k \mu_k \phi_k, \quad k_e \in {\mathbb N} .
\end{equation}

Using the same basis and $\lambda_k$'s, we next define on the span of $\{\phi_k\}$
\begin{equation}
\label{eq:Acc}
  A_{cc}\sum_{k=1}^{k_e} \alpha_k \phi_k = \sum_{k=1}^{k_e} \alpha_k \overline{\lambda_k} \phi_k.
\end{equation}
Hence $A_{cc}$ has the same eigenvectors as $A$, but its eigenvalues are the complex conjugates of those of $A$. 
Since (\ref{eq:var}) also holds for $\overline{\lambda_k}$, we have that 
\begin{equation}
\label{eq:eAcct}
  T_{cc}(t)\sum_{k=1}^{k_e} \alpha_k \phi_k = \sum_{k=1}^{k_e} \alpha_k e^{\overline{\lambda_k} t}\phi_k, \quad k_e \in {\mathbb N} 
\end{equation}
is the uniformly bounded semigroup generated by $A_{cc}$. Furthermore, for the same reason $A_{cc}$ is a bounded operator on $X$, with a densely defined inverse, $A_{cc}^{-1}$. Assume now that $A^{-1}$ and $A_{cc}^{-1}$ are both the infinitesimal generator of a $C_0$-semigroup. Then by standard semigroup theory, we have that $A_{cc}^{-1} +2I$ also generates a $C_0$-semigroup.
The expression of this generator is given by
\begin{align*}
  \left[A_{cc}^{-1} +2I\right]\sum_{k=1}^{k_e} \alpha_k \phi_k =&\ \sum_{k=1}^{k_e} (\overline{\lambda_k})^{-1} \alpha_k \phi_k  + \sum_{k=1}^{k_e} 2\alpha_k \phi_k \\
  =&\ \sum_{k=1}^{k_e} \overline{\mu_k} \alpha_k \phi_k  + \sum_{k=1}^{k_e} 2\alpha_k \phi_k = \sum_{k=1}^{k_e} (1 - i y_k) \alpha_k \phi_k\\
   =&\ 
  -A^{-1} \sum_{k=1}^{k_e} \alpha_k \phi_k .
\end{align*}
Hence if $A^{-1}$ and $A_{cc}^{-1}$ are both the infinitesimal generator of a $C_0$-semigroup, then $A^{-1}$ and $-A^{-1}$ are infinitesimal generators. This implies that $A^{-1}$ generates a $C_0$-group on $X$.

By \cite{Zwar10} we know that this cannot hold, since the eigenvalues of $A^{-1}$ are in a strip parallel to the imaginary axis and have uniform gap, and the eigenvectors do not form a unconditional basis.  So either $A^{-1}$ or $A_{cc}^{-1}$ (or both) cannot be  the infinitesimal generator of a $C_0$-semigroup, providing us 
the counter example.

We close this section with some remarks.
\begin{enumerate}
\item From the example it is very clear that the sequences of \cite{LoMV26}, repeated in (\ref{eq:y}) and (\ref{eq:lambda}) are the main ingredients. Since in the paper there is also a special conditional basis constructed, that might be less obvious from the examples in \cite{LoMV26}.
\item
  We know that $A^{-1}$ or $A_{cc}^{-1}$ does not generate a $C_0$-semigroup on $X$. For the sake of the argument, assume that it is $A^{-1}$. It is easy to see that its candidate semigroup is given on the span of the $\phi_k$'s as
  \begin{equation}
\label{eq:eA-1t}
  S(t)\sum_{k=1}^{k_e} \alpha_k \phi_k = \sum_{k=1}^{k_e} \alpha_k e^{\mu_k t}\phi_k, \quad k_e \in {\mathbb N} .
\end{equation}
If there would exists a $t_0>0$ such that $S(t)$ would be a (uniformly) bounded operator on the interval $(0,t_0)$, then it is easy to show that $\left(S(t)\right)_{t\geq 0}$ would be a $C_0$-semigroup. Hence there exists none such $t_0$.
\item
  Using the conditional basis of \cite{LoMV26}, we find that $S(\pi)$ is an unbounded operator. Since $S(\pi) = S(\pi/n)^n$, $n \in {\mathbb N}$, we have that for that basis $S(\pi/n)$ is an unbounded operator for all $n \in {\mathbb N}.$
\item
  I think that both $A^{-1}$ and $A_{cc}^{-1}$ do not generate a $C_0$-semigroup on $X$.
\end{enumerate}

\section{Some personal comments}

Looking at the above counter example, I wondered why I did not constructed it 10--20 years ago. Of course it is very hard to remember what I did and did not try to solve this problem. However, I have a clear idea why I completely overlooked this. The choice for a conditional basis is a natural one when you want to construct counter example, and I had gained quite some experience with it. Naturally, I also knew the result of \cite{Zwar10} stating that you cannot construct a $C_0$-group with spectral gap and a conditional basis. However, my approach was always to start with an unbounded generator of a bounded semigroup, and to show that its (bounded) inverse would not generate to a bounded semigroup. Hence starting from an operator that does not generate a $C_0$-semigroup never crossed my mind.  As is clear from the example in the previous section, this is not exactly what is done. The example starts on the ``bounded side'' and next study the unbounded inverse. Looking back, I don't know why I never tried starting with a generator of a ``non-$C_0$-semigroup'' and prove that its inverse generated a bounded semigroup. After all, by \cite{Zwar07} the problem of ``always bounded'' is equivalent to ``always infinitesimal generator''.

\section{More on the eigenvalues}

Similar as in (\ref{eq:lambda}) we take $\mu_k = -1+ i z_k$ and  $\lambda_k=\mu_k^{-1}$, but we study different sequences than $y_k$'s. We assume on the positive sequence $z_k$ the following; there exits an $\varepsilon >0$ such that $z_{k+1} > z_k +\varepsilon$ for all $k$. Thus the sequence $\mu_k$ has an uniform gap. 

If we want that $A$ as in (\ref{eq:A}) generates a (uniformly) bounded semigroup, then it is sufficient to show that
\begin{equation}
\label{eq:bnd-el}
  \sup_{t\geq 0} \sum_{k=1}^{\infty} \left| e^{\lambda_{k+1} t} -  e^{\lambda_{k} t}\right| < \infty.
\end{equation}

Next we study this infinite sum. We distinguish two cases the exponential $e^{\lambda_{k} t}$ close to 0 and the exponential $e^{\lambda_{k} t}$ close to 1. We fix $\alpha >0$ and for a fixed $t>0$ we define
\begin{equation}
\label{eq:k0}
  k_0 \in {\mathbb N} \mbox { such that } \mathrm{Re}(\lambda_{k_0}t) \leq -\alpha \mbox{ and } \mathrm{Re}(\lambda_{k_0+1}t) > -\alpha. 
\end{equation}
Note that 
\begin{equation}
\label{eq:Relamb}
   \mathrm{Re}(\lambda_{k}t) = \frac{-t}{1+z_k^2}
\end{equation}
and by the conditions on $z_k$ this is increasing and converging to zero over the negative real axis.

We spit the infinite sum in (\ref{eq:bnd-el}) at $k_0$, and estimate these partial sum differently.
\begin{equation}
\label{eq:sum1}
  \sum_{k=1}^{k_0} \left| e^{\lambda_{k+1} t} -  e^{\lambda_{k} t}\right| \leq 2 \sum_{k=1}^{k_0} \left| e^{\lambda_{k} t} \right| = 2 \sum_{k=1}^{k_0} e^{\mathrm{Re}(\lambda_{k}) t}.
\end{equation}
\begin{equation}
\label{eq:sum2}
  \sum_{k=k_0}^{\infty} \left| e^{\lambda_{k+1} t} -  e^{\lambda_{k} t}\right| \leq  \sum_{k=k_0}^{\infty} \left| e^{\lambda_{k} t} \right| |\lambda_{k+1} t -  \lambda_{k} t |\leq  \sum_{k=k_0}^{\infty}  |\lambda_{k+1} t -  \lambda_{k} t |.
\end{equation}

We study these sum for different sequences $z_k$. To make life easy, we take $k_0$ s.t.\ 
\[
  \frac{-t}{1+z_{k_0}^2} = -\alpha \Leftrightarrow 1+z_{k_0}^2 = t/\alpha.
\]
\subsection{For the sequence $z_k=k^p$, $p>1$}

For this sequence of $z_k$'s, we have that 
\[
  \lambda_{k+1}- \lambda_{k} = \frac{i(z_k-z_{k+1})}{(-1+i z_{k+1})(-1+i z_{k})}.
\]
Its norm behaves like $k^{-p-1}$.
From (\ref{eq:k0}) and (\ref{eq:Relamb}) we have that $k_0 \approx \sqrt[2p]{t/\alpha}$. With the integral criterion we find that (\ref{eq:sum2}) behaves like
\[
  \sum_{k=k_0}^{\infty}  |\lambda_{k+1} t -  \lambda_{k} t | \approx t \sum_{k=k_0}^{\infty}  k^{-p-1} \approx t k_0^{-p} \approx \sqrt{t}.
\]

\subsection{For the sequence $z_k = \sqrt{e}^k$}

For this sequence of $z_k$'s, we have that 
\[
  \lambda_{k+1}- \lambda_{k} = \frac{i(z_k-z_{k+1})}{(-1+i z_{k+1})(-1+i z_{k})}.
\]
Its norm behaves like $\sqrt{e}^{-k}$. From (\ref{eq:k0}) and (\ref{eq:Relamb}) we have that $k_0 \approx \log(t/\alpha)$. With the integral criterion we find that (\ref{eq:sum2}) behaves like
\[
  \sum_{k=k_0}^{\infty}  |\lambda_{k+1} t -  \lambda_{k} t | \approx t \sum_{k=k_0}^{\infty} \frac{1}{\sqrt{e}^{k} } \approx t \frac{1}{\sqrt{e}^{k_0} } = \sqrt{t}.
\]
We find that (\ref{eq:sum1}) behaves like
\[
  \sum_{k=1}^{k_0}  e^{\mathrm{Re}(\lambda_{k}) t} \approx \sum_{k=1}^{k_0}  e^{-e^{-k} t} 
\]
With the integral criterion and making the integral substitution $\xi = e^{-k}t$, and we find
\[
  \int_{k=1}^{k_0}  e^{-e^{-k} t} dk = \int_{\xi=t/e}^{\alpha}  e^{-\xi} \frac{-1}{\xi} d\xi
\]
{\bf Remark:} If, instead of $k_0$, we take the ``breaking index'' at $2k_0$, then (\ref{eq:sum2}) becomes bounded, but (\ref{eq:sum1}) behaves like $\log(t)$.

Based on the above, we try another sequence.

\subsection{For the sequence $z_{k+1} = z_k^p$ with $p>1$ and $z_1=e$.}

For this sequence of $z_k$'s, we have that 
\[
  \lambda_{k+1}- \lambda_{k} = \frac{i(z_k-z_{k+1})}{(-1+i z_{k+1})(-1+i z_{k})}.
\]
Its norm behaves like $z_k^{-1}$. We find that 
\[
  \log(z_k) = p^{k-1},
\]
and so from (\ref{eq:k0}) and (\ref{eq:Relamb}) we have that $k_0 \approx \log_p\left(\log(t/\alpha)\right)$. With the integral criterion we find that (\ref{eq:sum2}) behaves like
\[
  \sum_{k=k_0}^{\infty}  |\lambda_{k+1} t -  \lambda_{k} t | \approx t \sum_{k=k_0}^{\infty} \frac{1}{e^{p^{k-1}} } \approx t \frac{e^{-p^{k_0}}}{p^{k_0}} .
\]
This last expression is uniformly bounded. 
Next we study (\ref{eq:sum1}).
\[
  \sum_{k=1}^{k_0}  e^{\mathrm{Re}(\lambda_{k}) t} \approx \sum_{k=1}^{k_0}  e^{-\frac{z_{k_0}^2}{z_k^2}} 
\]
Again by the integral criterion, this is uniformly bounded. 

\subsection{Summary}
What we shown with the above analyses that if the sequence $z_k$'s grows fast enough, the operator $A$ of (\ref{eq:A}) generates a bounded semigroup. For a sequence like $z_k=a^k$ we could show a log-growth, but not uniformly bounded. Only for double exponential sequences, like was chosen in (\ref{eq:y}), we could show that the corresponding semigroup was uniformly bounded. 

\subsection*{Acknowledgement} The author wants to thank Emil Lorist and Mark Veraar very much for the useful discussions leading up to this manuscript.

\end{document}